
\magnification=\magstep1
\input epsf
\parindent=0in
\parskip=.1in

\def\e{\varepsilon}

\def\m{\mu}

\def\r{\rho}
\def\si{\sigma}

\def\Z{{\bf Z}}
\def\Q{{\bf Q}}
\def\R{{\bf R}}

\def\A{{\bf A}}
\def\P{{\bf P}}

\def\ZZ{{\bf Z}}
\def\QQ{{\bf Q}}
\def\RR{{\bf R}}
\def\CC{{\bf C}}

\def\PP{{\bf P}}

\def\i{\infty}

\def\p{\prod}
\def\s{\sum}

\def\sub{\subseteq}
\def\ra{\rightarrow}

\def\smallsetminus{\setminus}

\def\D{\Delta}

\def\cD{{\cal D}}

\def\cU{{\cal U}}

\def\cF{{\cal F}}

\def\cI{{\cal I}}

\def\v{\vskip .1in}

\def\[{$$}
\def\]{\[}

\def\pl{\partial}

\def\sets#1{{[\![#1]\!]}}

\def\Vol{\mathop{\rm Vol }\nolimits}

\def\Norm{{\rm N }}
\def\Trace{\mathop{\rm Tr }\nolimits}
\def\Star{\mathop{\rm St }\nolimits}
\def\Link{\mathop{\rm Lk }\nolimits}

\def\ker{\mathop{\rm ker }\nolimits}
\def\sign{\mathop{\rm sign }\nolimits}
\def\dim{\mathop{\rm dim }\nolimits}

\newcount\snum
\newcount\unum
\newcount\thmnum
\newcount\lemnum
\newcount\propnum
\newcount\defnum
\newcount\eqnum
\snum=0
\unum=0
\thmnum=0
\lemnum=0
\propnum=0
\defnum=0
\eqnum=0

\newcount\lnum
\lnum=1
\def\slist{ \smallskip\lnum=1\narrower\narrower\narrower\narrower}
\def\elist{ \smallskip\noindent}
\def\lister{\par\noindent { $({\the\lnum})$} \advance\lnum by 1}

\font\frteenrm=cmr10 at 20truept
\font\fsteenrm=cmr12
\font\abrm=cmr10

\def\today{\ifcase\month\or
  January\or February\or March\or April\or
  May\or June\or July\or August\or
  September\or October\or November\or December\fi
  \space\number\day, \number\year}

\def\section#1{\v\v\goodbreak\advance\snum by 1{\bf \fsteenrm {\the\snum}.\ \ #1}
\hfil \unum=0}
\def\subsection#1{\v \advance\unum by 1{\bf \the\snum.\the\unum.}\ }
\def\showeqnum{\global\advance\eqnum by 1{\the\eqnum}}
\def\eqref#1{{(#1)}}
\def\ref#1{}
\def\emph#1{{\it #1\/}}
\def\cite#1{{{\rm[#1]}}}

\def\theorem{\advance\thmnum by 1 {\bf Theorem {\the\thmnum}. }\it}
\def\endtheorem{\rm}

\def\lemma{\advance\lemnum by 1 {\bf Lemma {\the\lemnum}. }\it}
\def\endlemma{\rm}

\def\proposition{\advance\propnum by 1 {\bf Proposition {\the\propnum}. }\it}
\def\endproposition{\rm}

\def\definition{\advance\defnum by 1 {\bf Definition {\the\defnum}. }}
\def\enddefinition{\rm}

\def\qed{\vrule height4pt width4pt}
\def\proof{{\it Proof. }}
\def\endproof{\hfill\qed}

\def\innprod#1#2{\langle #1, #2 \rangle}
\def\frac#1#2{{#1 \over  #2}}

\def\title{Units, polyhedra, and a conjecture of Satake}

\nopagenumbers

\newif\iftitle
\def\titlepage{\global\titletrue}


\headline={\iftitle\titleheadline
                \else\pageheadline\fi}
\footline={\iftitle\global\titlefalse\titlefootline
                   \else\pagefootline\fi}

\def\titleheadline{\hfil}
\def\titlefootline{\hfil\folio\hfil}

\def\pagefootline{\hfil}
\def\pageheadline{\ifodd\pageno\rightheadline \else\leftheadline\fi}
\def\rightheadline{\tenrm\hfil{\it \title}\quad\quad\folio}
\def\leftheadline{\tenrm\folio\quad\quad Paul E. Gunnells and Jacob Sturm\hfil}
\voffset=2\baselineskip

\titlepage

\newif\ifcentered
\centeredfalse

\ifcentered
\bigskip
\bigskip
\centerline{{\frteenrm \title}}
\bigskip
\bigskip
\centerline{Paul E. Gunnells and Jacob Sturm}
\centerline{Dept. of Mathematics, Rutgers University, Newark NJ 07102--1811}
\bigskip
\centerline{\today}
\bigskip
\bigskip
\else
\bigskip\bigskip
\line{{\frteenrm \title}\hfill }
\bigskip\bigskip
\line{Paul E. Gunnells\footnote{}{Supported in part by the National
Science Foundation under grants DMS 00-70747 and DMS 01-00410.} and Jacob Sturm\hfill }
\medskip
May 21, 2002
\bigskip\bigskip\bigskip
\line{{\fsteenrm Abstract}\hfill }
\bigskip
{\abrm
Let $F/\QQ $ be a totally real number field of degree $n$.  We
explicitly evaluate a certain sum of rational functions over a
infinite fan of $F$-rational polyhedral cones in terms of the norm map
$\Norm \colon F\rightarrow \QQ $.  This completes Sczech's
combinatorial proof of Satake's conjecture connecting the special
values of $L$-series associated to cusp singularities with
intersection numbers of divisors in their toroidal resolutions
\cite{Sc2}.
}
\bigskip
\bigskip
\fi

\section{Introduction}
\subsection{}
Let $F/\QQ $ be a totally real number field of degree $n$ with ring of
integers $\ZZ _{F}$.  Let $M\subset F$ be a $\ZZ $-module of rank $n$,
and let $\rho \in F$ determine a coset $M+\rho $.  Recall that $x\in
F$ is said to be \emph{totally positive} if $v (x) > 0$ for each infinite
place $v$ of $F$.  Let $U\subset \ZZ _{F}^{\times }$ be the totally
positive units, and let $U_{M+\rho }\subset U$ be the subgroup with
$\varepsilon (M+\rho ) = M+\rho $.   Choose a finite index
subgroup $V\subset U_{M+\rho }$.  Then the pair $(M+\rho , V)$ determines a
collection of special values
\[
L (M+\rho ,V;s) := \sum _{\mu \in (M+\rho )/V} \Norm (\mu )^{-s},
\quad s=1,2,3,\dots .
\]
Here $\Norm \colon F\rightarrow \Q$ is the norm map, and in the sum we
omit $\mu  = 0$.  These series converge
absolutely if $s>1$, but only conditionally if $s=1$.  The theorem of
Klingen and Siegel asserts that these special values are---up to $d
(M)$, the square root of the discriminant of $M$, and powers of
$\pi$---cyclotomic numbers, and moreover rational if $\rho \in M$.
For example, let $F=\QQ (\sqrt{3})$, $M=\ZZ 1 + \ZZ (\sqrt{3}/3)$,
$\rho =0$, and $V=U_{M}=U$.  Then $d (M) = 2/\sqrt{3}$, and by
computing the $L$-series numerically to high accuracy (using GP-Pari
\cite{GP}) we easily find
\[
L (M,V;1) = -\frac{\pi ^{2}\sqrt{3}}{6}, \quad L (M,V;2) = \frac{\pi
^{4}\sqrt{3}}{6}, \quad L (M,V;3) = -\frac{\pi ^{6}\sqrt{3}}{36}.
\]

In \cite{Sa}, Satake proposed a geometric interpretation of these
special values in terms of certain intersection numbers.  Associated
to the data $(M,V,s)$ is a cusp singularity $Y^{*} (M,V,s)$; this is a
complex analytic space of complex dimension $ns$ with an isolated
singular point (the cusp), and is obtained by (partially)
compactifying the quotient of a bounded symmetric domain by an
arithmetic subgroup.  Examples of this class of singularities include
the cusps in the Baily-Borel compactifications of Hilbert and Picard
modular varieties.  Using toroidal geometry \cite{AMRT, KKMS}, one can
construct a desingularization $\tilde{Y} (M,V,s)\rightarrow Y^{*}
(M,V,s)$ with an exceptional divisor with finitely many irreducible
components $\{D_{\tau } \}$; for $s=1$ these components are toric
varieties, and for $s>1$ they are toric variety bundles over an $n
(s-1)$-dimensional abelian variety $A$.  This abelian variety depends
on additional data used to construct $ \tilde{Y}(M,V,s)$, but the
exact nature of this data is not important for us.  Then Satake's
conjecture, in the special case $\rho \in M$, becomes
\[
\frac{\Vol (A)  d (M) ((s-1)!)^{n}}{(2\pi i)^{ns}}
\,L (M,V;s) =
\frac{1}{(ns)!} \,\bigl(\sum_{\tau }B_{\tau }D_{\tau
}\bigr)^{ns}.\eqno{({\rm S})}
\]
Here $B_{\tau }$ is a ``Bernoulli symbol'' defined such that $B_{\tau
}^{k}$ is the $k$th Bernoulli number $B_{k}$, independent of $\tau $.
The modification to $\rho \not \in M$ is slightly more complicated to
state, but the geometry remains unchanged.  Continuing the example
from above, we find that the resolutions $\tilde{Y} (M,V;s)$ have
exceptional divisors $D_{s}$ with two irreducible components $D_{s,1},
D_{s,2}$.  Computing the relevant intersection numbers yields
$$
\kappa D_{1,1}^{2} = -2, \quad \kappa D_{1,2}^{2} = -3, \quad \kappa D_{1,1}D_{1,2} = 2,
$$
$$
\kappa D_{2,1}^{4} = 0, \quad \kappa D_{2,2}^{4} = 0, \quad \kappa D_{2,1}^{2}D_{2,2}^{2} = 3,
$$
$$
\kappa D_{3,1}^{6} = -12, \quad \kappa D_{3,2}^{6} = -81/2, \quad \kappa D_{3,1}^{4}D_{3,2}^{2}
= -18, \quad \kappa D_{3,1}^{2}D_{3,2}^{4} = -27,
$$
where $\kappa $ denotes the inverse volume of the appropriate abelian
variety.  With these intersection numbers, \eqref{S} becomes
\[
\eqalign{
L (M,V;1) &=-\frac{4\pi^{2}}{2!}\frac{\sqrt{3}}{2}\Bigl (
-3B_{2}+4 \bigl(B_{1})^{2} - 2B_{2}\Bigr) = -\frac{\pi^{2}\sqrt{3}}{6},\cr
L (M,V;2) &=\frac{16\pi^{4}}{4!}\frac{\sqrt{3}}{2}\Bigl (
0\cdot B_{4} + 6 \cdot 3\bigl(B_{2}\bigr)^{2} + 0\cdot B_{4}\Bigr) = \frac{\pi^{2}\sqrt{3}}{6},\cr
L (M,V;3) &=-\frac{64\pi^{6}}{6!(2!)^{2}}\frac{\sqrt{3}}{2}\Bigl (
-12B_{6} - 15\cdot18B_{4}B_{2} - 15\cdot 27
B_{2}B_{4}-\frac{81}{2}B_{6}\Bigr) = -\frac{\pi^{6}\sqrt{3}}{36},\cr
}
\]
where we use the Bernoulli numbers $B_{1}=-1/2$, $B_{2} =
1/6$, $B_{4}= -1/30$, $B_{6}=1/42$, and where we omit $B_{0}=1$.

\subsection{}
Proofs of this conjecture already appear in the literature in special
cases, and always for $\rho \in M$.  For $s=1$ and $F$ real quadratic,
the conjecture for the $L$-value $L (M,V;1)$ was proved by Hirzebruch
\cite{Hir}, as a by-product of his resolution of the cusp
singularities of Hilbert modular surfaces.  For all totally real
fields and $s=1$, independent proofs were given by
Atiyah-Donnelly-Singer \cite{ADS} and M\"uller \cite{M\"ul}.

The most general previous results are due to Ogata
\cite{Oga}, who verified \eqref{S} for all $n$ and odd $s$.  The
restriction to odd $s$ arises as follows.  Like \cite{ADS} and
\cite{M\"ul}, Ogata considers the $L$-function
\[
L^{**} (M,V;s) := \sum _{\mu \in
M/V} \sign (\Norm (\mu ) ) |\Norm (\mu ) |^{-s},
\]
which converges in a suitable halfplane, has an analytic continuation
to the entire complex plane, and satisfies a functional equation of
the form $s \leftrightarrow 1-s$.  He then studies special values of
$L^{**}$ at nonpositive integral $s$.  The special values of $L^{**}$
coincide with those of $L$ for odd positive integral $s$, but not for even
positive integral $s$---in fact for even positive $s$ the series
$L^{**}$ has a much more complicated special value (analogous to the
special values of the Riemann zeta function at odd integers $\geq 3$),
to which the theorem of Klingen and Siegel doesn't apply.

The proofs of \cite{Hir} and \cite{ADS, M\"ul, Oga} differ in
essential ways.  Hirzebruch's proof relied on previous work of Meyer
\cite{Mey, Zag}, who had computed the $L$-values $L (M,V;1)$ for real
quadratic $F$ explicitly in terms of Dedekind sums.  Hirzebruch was
able to compute (essentially) the right of \eqref{S} explicitly by a
very concrete topological argument, and observed that the two sides
matched up to a constant.  The proofs of Atiyah-Donnelly-Singer,
M\"uller, and Ogata, by contrast, use much heavier machinery: the
Atiyah-Patodi-Singer index theorem, applied to the the manifold with
boundary obtained by ``cutting off'' the cusp singularity.  The main
drawback to either to these approaches is that \eqref{S} appears to be
true by sheer magic, and one doesn't directly see why the intersection
numbers have anything to do with special values of $L (M,V;s)$.

\subsection{}
In \cite{Sc2} Sczech describes a proof of Satake's conjecture that
avoids the difficult analysis of the index theorem, and explicitly
shows how the special values are built from the intersection numbers
of the cusp divisors.  The first step of the proof is to perform a
``Hecke unfolding" of the $L$-series, from a sum over $(M+\r)/V$ to a
sum over the full coset $M+\r$. This is done using the following
identity:
$$ L(M,V;s) \ = \ \sum_{\m\in(M+\r)/V}{1\over N(\m)^s}\ = \
\sum_{\m\in (M+\r)/V}\sum_{\sigma \in C} f_s(\sigma, \m)
\eqno(\bullet)
$$
where $C$ is a certain $V$ invariant $F$-{\it rational polyhedral fan}
supported in $(\R_{>0})^n$, and
the function $f_s(\sigma,x)$ is, for each fixed positive
integer  $s$ and each $\sigma\in C$,
a certain rational function in $n$ variables which Sczech constructs
explicitly from the data defining the fan.

$V$-equivariance properties of the rational functions $f_{s}$ allow
one to write the outer sum in \eqref{$\bullet$} as a sum over $M+\rho$
and the inner sum as a sum over $C/V$.  This shows that the $L$-value
can be rewritten as a finite sum of ``full sums," where a full sum is
a sum over the affine lattice $M+\rho $ of certain rational functions.
Each full sum is equal, in turn, to an Eisenstein series in the style
of \cite{Sc1}. The proof of \eqref{S} then consists of a careful
evaluation of the Eisenstein series using the formulas derived in
\cite{Sc1}.

The main purpose of this paper is to provide a proof of
\eqref{$\bullet$}.  Together with the arguments in \cite{Sc2}, this
completes the proof of Satake's conjecture.  Actually we prove an
apparently weaker statement (Theorem 1): we consider only $s=1$ and
$x_{0}$ totally positive, and take the sum in \eqref{$\bullet$} only
over the top dimensional cones of $C$.  However, easy arguments show
that this suffices for the application in \cite{Sc2}.

\newcount\seclistnum
\seclistnum=1 
\def\bumpseclistnum{\advance\seclistnum by 1\S\the\seclistnum}
\def\showseclistnum{\S\the\seclistnum}
%

The convergence of \eqref{$\bullet$} is delicate.  Since the sum does
not have a natural ordering, it has a meaningful value only if it
converges absolutely. In \cite{Sc2} the absolute convergence of the
sum is established for all $x_0 \in F \sub \R^n$ (here we regard $F$
as a dense subset of $\R^n$ via its $n$ embeddings into $\R$). One
might naturally ask if the sum converges absolutely for all $x_0 \in
\R^n$, but one can show this is not the case: there is a dense
subset $\Sigma \sub \R^n$ such that the sum in \eqref{$\bullet$} fails
to converge absolutely for all $x_0 \in \Sigma$.

For the rest of this introduction we sketch the proof of Theorem 1 and
give a guide to the contents of our article.  For $s=1$ and $\sigma $
a top-dimensional cone $t$, we have $K=1$, and the rational function
$f_{s} (\sigma ,x)$ becomes a function that we denote $h^{*} (t)
(x)$.  Hence we must show
$$
\sum _{t} h^{*} (t) (x_{0}) = \frac{1}{\Norm(x_{0})},\eqno{(\bullet_{1})}
$$
where the sum is taken over all top-dimensional cones $t\in C$, and
where $x_{0}$ is totally positive.

The rational function $h^{*} (t) (x)$ satisfies a certain signed
additivity property with respect to simplicial subdivisions of
top-dimensional cones (Proposition 1), and $h^{*} (\Delta _{\infty })
(x_{0})$ is exactly equal to $1/\Norm (x_{0})$.  Hence a natural idea
to prove \eqref{$\bullet_{1}$} is to construct a sequence of partial
sums of rational functions $\{S_{N} (x) \}$ by constructing a finite
sequence of unions of top-dimensional cones $\{\Sigma _{N} \}$ such
that $\Sigma _{N} \rightarrow \Delta _{\infty }$, and then by defining
$S_{N} (x)$ by summing $h^{*}$ over the cones in $\Sigma _{N}$:
$$
S_{N} (x) = \sum _{t\in \Sigma _{N}} h^{*} (t) (x).
$$
One can then try to show that the sequence $S_{N} (x_{0})\rightarrow
h^{*} (\Delta _{\infty }) (x_{0})$ as $N\rightarrow \infty $ by
applying the additivity relation.

The main difficulty in carrying this out is that the value of the
partial sum $S_{N} (x_{0})$ depends rather subtlely on the geometric
properties of $\Sigma _{N}$.  For example, the singular hyperplanes of
the rational function $S_{N} (x)$ are the linear spans of the cones in
the boundary $\partial \Sigma _{N}$ of $\Sigma _{N}$, and so one must
construct $\Sigma _{N}$ so that these linear spans are far from
$x_{0}$ for large $N$.  Furthermore, if $\partial \Sigma _{N}$ is
complicated there is highly nontrivial cancellation among the terms in
$S_{N} (x_{0})$, which makes a direct estimate of the error term
infeasible.

To resolve these difficulties, we first define a notion of
\emph{cycles} built from rational polyhedral cones, and show that
$h^{*}$ induces a cocycle with coefficients in the field of rational
functions in $n$ variables (\S\S 3--7).  This allows us to work with
``formal polyhedra'' instead of actual convex polyhedra, and to avoid
explicitly dealing with cancellations.  Then we analyze the action of
the unit group $V$ on $\Delta _{\infty }$ to construct a sequence of
convex polyhedral cones $\{\Sigma _{N} \}$ exhausting $\Delta _{\infty
}$, and such that the singular hyperplanes spanned by $\partial \Sigma
_{N}$ are far from $x_{0}$ (\S 8).  Finally we apply a geometric
interpretation of $h^{*}$ due to Hurwitz \cite{Hur} to complete the
proof of Theorem 1 (\S 9).

\subsection{Acknowledgements.}
We thank Robert Sczech for suggesting this problem to us, for giving
us access to \cite{Sc2}, and for many interesting conversations.  We
thank the NSF for support, and the first named author is grateful for
the support of the Max Planck Institut in Bonn.

\section{Notation and statement of the main result}
\subsection{}
We retain the notation from \S 1.  In particular, $F/\QQ $ is a
totally real number field of degree $n$ with ring of integers $\ZZ
_{F}$, $M\subset F$ is a rank $n$ $\ZZ $-module, $\rho \in F$ gives a
coset $M+\rho $, and $V\subset U_{M+\rho } \subset U$ is a finite-index
subgroup of the totally positive units preserving $M+\rho $.

We fix an ordering $v_{1}, \dots , v_{n}$ of the infinite places of
$F$, and identify $F$ with its image in $\RR ^{n}$ given by these
embeddings.  We use the notation $x^{(i)}$ for $v_{i} (x)$.  Let
$\innprod{\phantom{a}}{\phantom{a}}\colon (\RR ^{n})^{2}\rightarrow
\RR $ be the standard scalar product.  Note that if $x,y\in F$, then
$\innprod{x}{y} = \Trace (xy)$.  We will use this scalar product to
identify $\RR ^{n}$ with its dual.

\subsection{}
We recall some notations about cones and fans.  For more information
we refer to \cite{Ful}.  A subset $\sigma $ of a real vector space $W$
is called a \emph{cone} if
\slist
\lister $\sigma \cap -\sigma = \{0 \}$, and
\lister $x\in \sigma $ implies $\lambda x\in \sigma $ for all $\lambda
\in \RR _{\geq 0}$.
\elist
For any cone $\sigma $, the \emph{dual
cone} $\sigma ^{*}$ is defined by
$$
\sigma ^{*} := \bigl\{ x \bigm|  \innprod{x}{y} \geq 0 \quad \hbox{for all $y\in
\sigma $} \bigr\}.
$$
A cone is \emph{polyhedral } if it
is the convex hull of finitely many half-lines, \emph{simplicial} if
the number of these lines can be taken to be the dimension of the
linear span $\RR \sigma $ of $\sigma $, and \emph{rational} if $W$ has a $\QQ
$-structure and these half-lines each contain nonzero points in $W (\QQ )$.
Any rational polyhedral cone $\sigma $ contains a collection of
subcones, the \emph{faces} of $\sigma $.  We write $\tau < \sigma $ to
indicate that $\tau $ is a face of $\sigma $.  Faces of codimension 1
are called \emph{facets}.

\subsection{}
Let $C$ be a set of rational polyhedral cones.  Then $C$ is a
\emph{fan} if
\slist
\lister $\sigma \in C$ and $\tau < \sigma $ implies $\tau \in C$, and
\lister if $\sigma ,\tau \in C$, then $\sigma \cap \tau $ is a face of
each.
\elist
We write $C (k)$ for the subset of $C$ consisting of $k$-dimensional
cones.  Following \cite{Sc2}, we reserve the notation $t$
(respectively $\tau $) for an element of $C (n)$ (resp. $C (1)$).

Given a cone $\sigma \in C$, the \emph{star} of $\sigma $ is the set
\[
\Star (\sigma ) := \bigl\{ \sigma ' \in C \bigm|  \sigma < \sigma ' \bigr\},
\]
and the \emph{link} of $\sigma $ is the set
\[
\Link (\sigma ) := \bigl\{ \sigma' \bigm|  \sigma' < t \quad \hbox{for
some $t\in \Star (\sigma )$, and $\sigma \not < \sigma '$}  \bigr\}.
\]

A fan $C'$ is a \emph{refinement} of $C$ if every cone in $C$ can be
written as a union of cones in $C'$.  Given any two rational
polyhedral fans, there exists a rational fan that is
a common refinement of each.

Recall that $\Delta _{\infty }$ is the totally positive chamber $(\RR_{>0})
^{n} \subset \RR ^{n}$.
Throughout this paper, we will only consider polyhedral
fans $C$ satisfying the following properties:
\slist
\lister $C$ is locally finite, apart from $\{0 \} \in C$,
\lister $C$ is simplicial,
\lister $C$ is $F$-rational, i.e. each $\sigma \in C$ is generated
by half-lines defined over $F$,
\lister there is a finite-index subgroup $V$ of the totally positive
units that
acts on $C$ via the embedding of $F$ into $\R^{n}$, and
\lister $C$ gives a decomposition of the totally positive chamber:
$$
\bigcup _{\sigma \in C \atop \sigma \not = \{0 \}} \sigma  = \Delta _{\infty }.
$$
Note that this last condition implies $|C|$ is infinite.  We call fans
satisfying conditions (1)--(5) \emph{good fans}.  The existence of
good fans was proved in great generality by Ash in Chapter II of
\cite{AMRT} (there, a good fan is called a \emph{$\Gamma$-admissible decomposition}).

\subsection{}
In what follows, we use the notational convention of multivariables.
In particular, let $x$ be a real multivariable with components
$x_{1},\dots ,x_{n}$, and let $A = ( A_{1},\dots ,A_{n})$ be an
$n$-tuple of vectors in $\RR ^{n}$.  Let $\CC (x)$ be the field of
rational functions in the variables $x_{1},\dots ,x_{n}$.  Our goal is
to define two maps
$$
\eqalign{
h, h^{*}\colon (\RR ^{n})^{n} &\longrightarrow \CC (x)\cr
                             A&\longmapsto h (A) (x), h^{*} (A) (x)\cr
}
$$
that will play an important role in the sequel.

First, we define the function $h (A) (x)\in \CC (x)$ by
$$
h (A) (x) = \frac{\det A}{\innprod{x}{A_{1}}\cdots \innprod{x}{A_{n}}}.
$$
Next we define the ``dual function'' $h^{*} (A) (x)$ as follows.
If
the set $A$ is linearly dependent we put $h^{*} (A) (x) = 0$.
Otherwise we let $B = ( B_{1}, \dots , B_{n})$ be the dual basis to
$A$ with respect to our inner product, and put
$$
h^{*} (A) (x) := h (B) (x).
$$

As a function of the $A_{i}$, both $h$ and $h^{*}$ are homogeneous of
degree $0$.  Moreover, if $E=(e_{1},\dots ,e_{n})$ is the canonical
basis of $\RR ^{n}$, then $h (E) (x) = h^{*} (E) (x) = 1/\Norm (x)$,
where for $x\in \RR ^{n}$ we put $\Norm (x) = \prod x_{i}$ (this
agrees with the usual norm map on $F\subset \R ^{n}$).  Each of these
function also enjoys a cocycle property:

\proposition
\cite{Sc2}
Given $n+1$ vectors $A_{0},\dots ,A_{n}$, let $A^{(i)}$ be the tuple
\[
(A_{0}, \dots , \hat A_{i}, \dots , A_{n}),
\]
where the hat means to
delete the $i$th component.  Then we have
$$
\sum _{i=0}^{n} (-1)^{i} h (A^{(i)}) (x) = \sum _{i=0}^{n} (-1)^{i}
h^{*} (A^{(i)}) (x) = 0.
$$
\endproposition

\subsection{}
We recall the geometric interpretation of $h$ due to Hurwitz \cite{Hur}.
Suppose $\innprod{x}{A_{i}} \not =0$ for $A_{1}, \dots , A_{n}$, and
$\det A \not =0$.  The
points $A_{i}$ determine $2^{n}$ simplicial cones in $\RR ^{n}$, and
hence $2^{n-1}$ regions in $\PP ^{n-1} = \PP (\RR ^{n})$.  Among these
regions is a unique region $R$ that misses the hyperplane $\{y\mid
\innprod{x}{y} = 0  \}$.  Then up to sign, the value of $h (A) (x)$ is
given by the integral
\[
\int _{R} \Omega _{x}, \quad \hbox{where} \quad \Omega _{x} =
\frac{(n-1)!}{\innprod{x}{y}^{n}}\sum _{i=1}^{n}
(-1)^{i-1}y_{i}\,dy_{1}\cdots \hat dy_{i} \cdots dy_{n}.
\]
Here $\int _{R} \Omega _{x}$ is equal to the euclidean area of the
region $R$ viewed as a subset of $\R^n$ as follows: $R \sub
\P^{n-1}\backslash x^\perp =\R^n$.
Moreover, the sign is determined by fixing an orientation on an affine
chart containing $R$.

\subsection{}
We extend the notation for $h$ and $h^{*}$ as follows.  Recall that $M
\subset \RR ^{n}$ is a lattice.  Given any $1$-cone $\tau \in C (1)$,
let $A_{\tau }\in \tau \cap M$ be the nonzero point closest to the
origin.  Fix a global orientation on $\RR ^{n}$.  Then if $t\in C
(n)$, we let $A_{t}$ be the $n$-tuple $(A_{\tau } \mid \tau <t)$,
where the points are positively ordered with respect to the
orientation.  Then we define
$$
h (t) (x) := h (A _{t}) (x),
$$
and define $h^{*} (t) (x)$ by using the dual basis to $A_{t}$.  Note
that the singular hyperplanes of $h^{*} (t) (x)$ are exactly the
linear spans of the facets of $t$.

\subsection{}
Now consider the infinite sum of rational functions
$$
S (C) (x) := \sum _{t\in C (n)} h^{*} (t) (x).
$$
We want to evaluate this sum at $x=x_{0} \in \Delta_{\infty
}$, but we must be careful because singular hyperplanes of some terms
may pass through $x_{0}$.  Let $\cU_{x_{0}} $ be the set of singular
cones determined by $x_{0}$:
\[
\cU _{x_{0}} := \bigl\{\sigma \in C \bigm| \dim(\sigma)<n,
x_{0}\in \RR \sigma , x_{0}\not \in \RR \sigma
'\quad \hbox{for all}\quad \sigma ' < \sigma  \bigr\}.
\]
If $\cU_{x_{0}} = \emptyset $, then every term in $S (C) (x_{0})$ is
well defined, and we say that $x_{0}$ is \emph{nonsingular } with
respect to $C$.  One can show in this case that if $x_{0}\in F$ then
$S (C) (x_{0})$ is absolutely convergent \cite{Sc2}.

On the other hand, if $\cU_{x_{0}} \not =\emptyset $, then we must
make the following modifications to $S (C)$.  For any $\sigma \in \cU
_{x_{0}}$, in $S (C) (x)$ we replace the sum of set of rational functions
\[
H_{\sigma } = \bigl\{h^{*} (t) (x) \bigm|  t \in
\Star (\sigma ) \cap C (n)\bigr\}
\]
with the finite sum $\Theta _{\sigma } := \sum _{h\in H_{\sigma }} h$.
We do this for all $\sigma \in \cU_{x_{0}} $, and add together the
resulting rational functions.

\proposition
All terms in the sum constructed by the above procedure are well-defined.
\endproposition

\proof Let $\sigma \in \cU _{x_{0}}$.  We must show that $x_{0}$ does
not lie in the singular hyperplanes of the function $\Theta _{\sigma
}$.  By Proposition 1, the singular hyperplanes of $\Theta _{\sigma }$
are the linear spans of the dimension $n-1$ cones in $\Link (\sigma)
$.  Thus it suffices to show that $\Link (\sigma )\cap \cU_{x_{0}} =
\emptyset $.

Suppose not and let $\sigma '$ be a member of this intersection.
By assumption,
$x_0 \in \R\sigma' \cap \R\sigma = \R(\sigma'\cap \sigma)$.
Since $\sigma ' \not > \sigma$, the minimality of $\sigma$ is contradicted.
\endproof

One can also show in this case that if $x_{0}\in F$, the
resulting series is again absolutely convergent \cite{Sc2}.

\lemma
Let $x_{0}\in F$, and let $C$ be a good fan.  If $C'$ is any other
good fan, then $S (C) (x_{0}) = S (C') (x_{0})$.
\endlemma

\proof First suppose that $C'$ is a refinement of $C$, and that
$x_{0}$ is nonsingular with respect to both $C$ and $C'$.  By the
cocycle property, each term in $S (C) (x_{0})$ is a finite sum of
terms from $S (C') (x_{0})$.  Since both sums converge absolutely, this
implies they have the same sum.  The singular case is handled by
considering sums over stars as above, and we leave the details to the
reader.

Now if $C'$ is any good fan, we can construct an $F$-rational common
refinement $C''$ of $C$ and $C'$.  This will also be a good fan,
although perhaps with respect to a different subgroup $V''$.
The previous argument shows $S (C) (x_{0}) = S (C'') (x_{0}) = S (C')
(x_{0})$, which completes the proof.
\endproof

We are now ready to state our main result:

\theorem Let $C$ be a good fan (\S 2.3), and suppose that $x_{0}\in
F$ is totally positive.  Then
$$
S (C) (x_{0}) = \frac{1}{\Norm (x_{0})}.
$$
\endtheorem

The proof will occupy the rest of the article.


\section{Chains and cycles in linear varieties}

\subsection{}
Let $X$ be a linear variety of dimension $n$, that is $X $ is an
algebraic variety biregular with $\P^n$.  Since the automorphism group
of $\P^n$ is $PGL(n+1)$, any concept that involves the linear
structure of $\P^n$ is meaningful for $X$. In particular, if we fix an
isomorphism $\phi: \P^n \ra X$, we may define a line in $X$ to be the
image under the map $\phi$ of a line in $\P^n$.  The definition does
not depend on the choice of $\phi$ and thus the concept of a line
(plane, hyperplane, etc.)  is intrinsic to $X$.

If $X$ is a linear variety, then we define $X^*$ to be the set of
hyperplanes in $X$. Then $X^*$ has a canonical structure of linear
variety.

\subsection{}
We wish to define $C_k(X)$ and $Z_k(X)$, the chains and cycles of
dimension $k$.

Let $C_{0} (X) = \ZZ [X]$, that is, $C_{0} (X)$ is the free abelian
group generated by the points of $X$.  Thus $c_0 \in C_0$ is a map
$c_0\colon X \ra \Z$ with the property $c_0(x)=0$ for all but finitely
many $x$.  Let $\pl_0 \colon C_0(X) \ra \Z$ be the map $\pl_0 (c_0(x))
= \s_{x\in X} c_0(x)$, and let $Z_0 = \ker \pl_0$.

We define $C_k(X)$ inductively:
$$
 C_k(X) \ = \ \bigoplus_{L^{(k)} \sub X} Z_{k-1}(L^{(k)}),
$$
where $L^{(k)}$ ranges over all linear subvarieties of dimension
$k$, and $ Z_{k-1}(L^{(k)})$ is the group of $k-1$ cycles in $L^{(k)}$.
A typical element of $C_k(X)$ has the form $c_k = c_k(X) =
(z_{k-1}(L^{(k)}))$, where $z_{k-1}(L^{(k)})$ is a $k-1$ cycle in
$L^{(k)}$ that is zero for all but finitely many $L^{(k)}$.

If $Y, X$ are linear varieties with $Y \sub X$, then we have a
canonical map $ Z_{k-1}(Y) \hookrightarrow Z_{k-1}(X)$, which we
denote $z_{k-1}(Y) \mapsto z_{k-1}(Y,X)$.  Using this we define the
differential $\pl_k\colon C_k(X) \ra C_{k-1}(X)$ by
$$
\pl_k(c_k) \ = \ \pl_k((z_{k-1}(L^{(k)})) \ = \
  \s_{L^{(k)}\sub X} \  z_{k-1}(L^{(k)},X),
$$
and put $Z_k(X) = \ker(\pl_k)$.  Note that $\pl_k\circ \pl_{k-1} = 0$.

\subsection{}
We now give some examples of our construction.

\slist
\lister The cycle group $Z_0(L)$ is generated by elements of the form
$w-v$ where $v,w \in L$ (here $L \sub X$ is any
linear subvariety). The only relations are
$(w-v) + (v-u) = (w-u)$ and $(w-v) = -(v-w)$,
where $u,v,w \in L$.
\lister The chain group $C_1(X)$ is generated by symbols of the
form $\langle w-v\rangle$, with $w,v \in X$. The only
relations are $\langle w-v\rangle = -\langle v-w\rangle$ and
$\langle w-v\rangle + \langle v-u\rangle = \langle w-u\rangle$ where $w,v,u \in X$
are collinear. If $c_1 = \s a_i\langle w_i-v_i\rangle\  \in \ C_1(X)$,
note that
$\pl_1(c_1) = \s a_i(w_i-v_i) \in C_0(X)$.

\lister Every element of $Z_1(X)$ can be written
in the form $\s_{i=1}^r \langle v_{i+1}-v_i\rangle$, where
$v_{r+1} = v_1$.
\elist

\subsection{}
Now we discuss simplicial cycles.
Let $ A_0,\dots, A_k \in X$ be linearly independent
(the concept of linear independence is intrinsic to $X$).
Then we define $ \si_{k-1}(A_0,\dots,A_k)(X) \in Z_{k-1}(X)$
as follows:

\slist
\lister If $k=1$, then $\si_0(A_0,A_1) := A_0-A_1$.
\lister If $k>1$, then
$$
\si_{k-1}(A)(X) =  \si_{k-1}(A_0,\dots,A_k)(X)
 :=  \s_{r=0}^k (-1)^r\si_{k-2}(A_0,\dots ,\hat A_r , \dots A_k)(L_r),
$$
where $L_r$ is the linear space spanned by $A_0,\dots,\hat A_r , \dots, A_k$.
\lister If $ A_0,\dots, A_k \in X$ are linearly dependent. then
$ \si_{k-1}(A_0,\dots,A_k)(X) := 0$.

\proposition

\slist
\lister $\si_{k-1}(A)(X) \in Z_{k-1}(X)$.
\lister For every permutation $\pi $ on $k+1$ letters, we have
$$
\si_{k-1}(A_{\pi(0)}, \dots , A_{\pi(k)})
\ = \ \sign(\pi) \si_{k-1}(A_0,\dots,A_k).
$$
\lister $Z_{k-1}(X)$ is spanned by simplices.

\endproposition

\proof
The proofs of (1) and (2) are standard. For (3), let
$z_{k-1}(X) = (z_{k-2}(L)) \in Z_{k-1}(X)$. Then  by
induction
$$ z_{k-2}(L) \ = \ \s_{A \in \A (L)}\si_{k-2}(A)(L),
$$
where $\A(L)$ is a finite set of $k$ tuples in $L$.
Now let $A_0 \in X$ be arbitrary.
Then we can check that
$$ z_{k-1}(X) \ = \ \s_L\s_{A \in \A (L)} \si(A_0,A)(X).
$$
\endproof

\section{CPD functions}

\subsection{}
Let $G$ be a torsion free abelian group and let
$f\colon  X^{k+1} \ra G$ be a function. Then we say that
$f$ is a \emph{CPD function} if it satisfies the \emph{cocycle}
property, the \emph{permutation} property, and the \emph{degeneracy}
property:

{\bf C.} For every $A_0,\dots, A_{k+1} \in X$, we have
 $$
\sum_{r=0}^{k+1} (-1)^r f(A_0,\dots, \hat A_r, \dots, A_{k+1}) = 0.
$$
{\bf P.} For every permutation $\pi $ on $k+1$ letters, we have
$$  f(A_{\pi(0)}, \dots , A_{\pi(k)}) = \sign(\pi) f(A_0,\dots,A_{k}).
$$
{\bf D.} If $(A_0,\dots,A_{k})$ are linearly dependent,
then $f(A_0,\dots,A_{k})  = 0 $.

\subsection{}
Here are some examples of CPD functions.
\slist
\lister
Let $X=\P^{n-1}$ and let $G = \R(x)$, the field of rational functions
in $n$ real variables, viewed as an abelian group with respect to addition.
The function $h \colon (\R ^{n})^{n} \rightarrow \R(x)$ from \S
2.4 is homogeneous in the points $(A_{1},\dots ,A_{n})$.  Hence it
induces a function $h \colon (\P ^{n-1})^{n}\rightarrow \R (x)$, and this is a CPD function.
\lister
Let $X$ be a linear variety of dimension $n$, let
$X^*$ be the dual variety, and
let $G= Z_{n-1}(X^*)$. We define a CPD function $D\colon  X^{n+1} \ra G$
as follows.  If $A_0,\dots,A_n$ are linearly dependent,
then $D(A_0,\dots,A_n) := 0$. Otherwise we let
$D(A_0,\dots,A_n)  := \si_{n-1}(B_0,\dots,B_n)$,
where $B_i\in X^*$ is the linear space spanned by
$A_0,\dots, \hat A_i, \dots, A_n$.
\elist

The following proposition shows that a $G$-valued CPD function $f$
extends to give a $G$-valued cocycle $\tilde{f}$ on the cycle group
$Z_{k-1} (X)$.  This will play a key role in the proof of Theorem 1.

\proposition Let $f\colon X^{k+1} \ra G$ be a CPD function.  Then
there exists a unique homomorphism $\tilde f\colon Z_{k-1}(X) \ra G$
satisfying $\tilde f(\si(A_0,\dots,A_k)(X)) = f(A_0,\dots,A_k)$ for
all $A_0,\dots,A_k \in X$.  \endproposition

\proof
Let $\tilde{f}$ be the function defined on simplices in $Z_{k-1} (X)$
as in the statement of the proposition.
Then $\tilde{f}$ can be defined on any
cycle $z\in Z_{k-1} (X)$ by writing $z=\sum_{A} \sigma (A) $, where the
$\sigma (A)$ are simplices, and then putting $\tilde{f} (z) = \sum
_{A} \tilde{f} (\sigma (A))$.  That this extension is well-defined
follows immediately from the cocycle property of the CPD function $f$.
\endproof

\definition
Suppose $D$ is the CPD function in the second example above.
Then the extension $\tilde D$ gives a map $Z_{k-1}(X) \ra
Z_{k-1}(X^*)$.  We shall write $\tilde D((z_{k-1})(X)) =
[z_{k-1}(X)]^*$, and shall say that $[z_{k-1}(X)]^*$ is the \emph{dual
cycle} to $z_{k-1}(X)$.
\enddefinition

\section{Convex Polyhedra}

\subsection{}
We now extend some of the basic notions of convex geometry in affine
space to the linear variety $X$.  An \emph{open half space} of $X$ is
a connected component of the set $X\setminus (H_1\cup H_2)$, where
$H_1,H_2 \sub X$ are two different hyperplanes.  A \emph{closed half
space } is the closure of an open half space.  A \emph{convex polyhedron}
$K\subset X$ is a finite intersection of closed half spaces with the
property $K\cap H = \emptyset$ for some hyperplane $H$.  Thus a
convex polyhedron is just the image of a usual compact convex polyhedron in
$\R^{n+1}\setminus \{0 \}$ via the map $\R^{n+1}\setminus \{0
\}\rightarrow \PP ^{n}\rightarrow X$.  Moreover, there is similarly an
obvious bijection between polyhedral cones in $\RR ^{n+1}$ and
convex polyhedra in $X$.

Let $H\sub X$ be a hyperplane.  If $v_1,\dots,v_r$ are points in the
affine space $\in X\setminus H$, then the convex hull of
$\{v_1,\dots,v_r\}$ in $X\setminus H$ (in the usual sense of convex
geometry) is a convex polyhedron in $X$.  We can also define the convex hull
of a set of points in $X$ without passing to an affine subset, simply
by applying the usual definition with our notion of half space.  If $r
= k+1$ and if the $v_i$ are linearly independent, then we call the
resulting convex hull a \emph{simplicial polyhedron}.

\subsection{}
Now let $K \sub X$ be a convex polyhedron with $\dim(K) = \dim(X) = k$.  Let
$K^* \sub X^*$ be the set of hyperplanes in $X$ that do not intersect
the interior of $K$. Then $K^*$ is called the \emph{dual} of $K$, and
is itself a convex polyhedron in $X^*$.  It is easy to check that if $K$ is a
convex polyhedron corresponding to the polyhedral cone $\sigma _{K}\in \RR
^{n+1}$, then the dual convex polyhedron $K^{*}$ corresponds to the dual cone
$\sigma _{K}^{*}$.

An equivalent definition of $K^*$ is as follows.  Let $h$ by any fixed
point in the interior of $K$. Then $h^*$, the set of hyperplanes of
$X$ that pass through $h$, is a hyperplane in $X^*$.  Let
$\Phi_1,\dots,\Phi_t$ be the faces of $K$. Then the $\Phi_i$ span
hyperplanes $H_i$ that do not contain $h$ and thus determine points
$H_i^*\in X^*\setminus h^*$.  Then $K^*$ is the convex hull of $F_i^*$
inside $X^* \setminus h^*$.

The faces of $K^*$ are convex polyhedra in 1--1 correspondence with the
non-degenerate vertices of $K$ (these are the vertices not contained
in any of the planes spanned by the proper faces of $K$). If $v$ is
such a vertex, then the face $F_v^*$ corresponding to $v$ is the
convex polyhedron in $v^*$ whose vertices are the hyperplanes spanned by the
faces of $K$ containing $v$.  Thus $F_v^*$ is the convex hull of those
vertices, taken inside $v^*\setminus h^*$.

\subsection{}
There is another procedure for constructing $F_v^*\sub v^*$ that is
useful for inductive proofs.  Let $v \in K$ be a
non-degenerate vertex.  Then there exists a hyperplane $H_v \sub X$,
called a \emph{vertex hyperplane}, with the following properties:

\slist
\lister $v \notin H_v$.

\lister If $F_v = K\cap H_v$, then $K\setminus H_v$ consists
of two connected components. One of those components has, as its
closure, the convex polyhedron given by the convex hull of $F_v$,
and the point $v$.
\elist
The convex polyhedron  $ F_v^*$ is the dual of $F_v$ in the sense of \S 5.2.
To be precise, we have a
canonical map $H_v^* \ra v^*$ that associates to a
hyperplane $L\subset H_v$ the hyperplane of $X$
spanned by $v$ and $L$. Then the dual of
$F_v$, which is a convex polyhedron in $H_v^*$, is identified with $F_v^*$
via the canonical map $H_v^* \ra v^*$.

\section{The cycle associated to a polyhedron}

\subsection{}
Let $K\sub X$ be an oriented convex polyhedron of dimension $k$.  Then we
associate to the boundary of $K$ a cycle $z_{k-1}(K)(X) \in
Z_{k-1}(X)$ as follows.  If $K$ is 1 dimensional, then $\pl(K)$ is
just a pair of points $x,y \in X$. The orientation allows us to assign
one of the points, say $x$, the value $+1$, and the other point $y$
the value $-1$. Then we define $z_0(K)$ to be $x-y$.

In general, we define $z_{k-1}(K)$ inductively as follows.  For each
face $F \sub K$, let $L_F$ be the linear space spanned by $F$. Then
$F$ is a convex polyhedron of dimension $k-1$ embedded in $L_F$. Thus
$z_{k-2}(F)(L_F)$ has been defined as a $k-2$ cycle in $L_F$. We
now define
$$z_{k-1}(K) = \s_F z_{k-2}(F)( L_F).$$

In particular,
$$  z_{k-1}(K^*) \ = \ \s_v z_{k-2}(F_v^*) \ = \
   \s_v z_{k-2}([K\cap H_v])^*).
\eqno(\showeqnum)
$$
A cycle of the form $z_{k-1}(K) $ is called a \emph{polyhedral cycle}.

At this stage we have two notions of duality for a polyhedral cycle $z
(K)$: the cycle $z (K^{*})$ associated to the dual convex polyhedron $K^{*}$,
and the dual cycle $[z (K)]^{*}$ constructed using the CPD function
$\tilde{D}$ from Example 2 of \S 4.2.  The following theorem asserts
that these two notions coincide.

\theorem
Let $K \sub X$ be a convex polyhedron. Then we have
$$  z_{k-1}(K^*) \ = \ [z_{k-1}(K)]^*.
$$
\endtheorem

\proof Without loss of generality, we may assume the faces of $K$ are
simplicial polyhedra.  Choose a point $u$ in the interior of $K$.
Then $K = \cup_{F} c(u,F)$ is a decomposition into simplicial
polyhedra, where $c(u,F)$ is the cone on $F$ with vertex $u$, and $F$
ranges over all faces of $K$.  We obtain
$$z_{k-1}(K) = \s_F z_{k-1}(c(u,F)),
$$
and thus
$$[z_{k-1}(K)]^* = \s_F [z_{k-1}( c(u,F))]^*.$$

On the other hand, we have
$$
z_{k-1}(K^*) \ = \ \s_v z_{k-2}([K\cap H_v])^*)
$$
Hence we must show
$$ \s_v z_{k-2}([K\cap H_v])^*) \ = \ \s_F [z_{k-1}( c(u,F))]^* \ = \
\s_F z_{k-1}( c(u,F)^*), \eqno(\showeqnum)
$$
where $H_v$ is a vertex hyperplane.  Here the second equality
follows since $c(u,F)$ is a simplicial polyhedron, and
our theorem is trivially true for such polyhedra.

Each side of (2) is a cycle in $X^*$, and thus each side has a $v^*$
component for all $v \in X$. We shall compare the $v^*$ components for
each $v$ and show they are equal.  Fix $v \in X$. We may
assume $v$ is a non-degenerate vertex of $K$, for otherwise the $v^*$
component of each side of (2) vanishes.
Applying (1) with $K$ replaced by $c(u,F)$, and substituting
the result in (2), we are reduced to showing
$$z_{k-2}([K\cap H_v]^*) \ = \
\s_{\{F\mid v \in F\} }z_{k-2}([c(u,F)\cap H_v]^*).$$
But we clearly have

$$z_{k-2}([K\cap H_v]) \ = \
\s_{\{F\mid v \in F\} }z_{k-2}([c(u,F)\cap H_v]),\eqno{(\showeqnum)}
$$
since $K\cap H_v = \cup_{\{F \mid v \in F\} } [c(u,F)\cap H_v]$ is
a decomposition into simplicial polyhedra.
Taking duals of both sides of (3), and using induction on $k$,
completes the proof.
\endproof


\section{Cycles and the function $h^{*}$}

\subsection{}
Now we want to apply the machinery in \S\S 3--6 to
compute the sum $S (C) (x_{0})$ for totally positive $x_{0}\in F$.
Let $X$ be the projective space $\PP ^{n-1}$.  Recall that we have
identified $\RR ^{n}$ with its dual.  This allows us to identify $X$
with its dual $X^{*}$, and thus to identify $Z_{k} (X)$ with $Z_{k}
(X^{*})$.  We will do this throughout the following discussion.  Also,
from this point on we will only need to consider the cycle group
$Z_{n-2} (X)$, and so to lighten notation we will drop the degree
subscripts if no confusion is possible.

\subsection{}
Let $t\in C (n)$ be a top dimensional cone.  The results from the
previous sections show that $t$ determines a polyhedron $K_{t}$ and a
cycle $z (K_{t}) \in  Z (X)$.  If $t^{*}$ is the dual cone with associated
polyhedron $K_{t}^{*}$, then we have
\[
(z (K_{t}))^{*} = z (K_{t}^{*}),
\]
where the star on the left denotes the dual cycle construction.
Let $\tilde{h}\colon Z (X) \rightarrow \CC (x)$ be the extension to $Z
(X)$ of the function
described in Example (1) from \S 4.2 (whose existence
is guaranteed by Proposition 4).  Then we have
\[
h^{*} (t) (x) = \tilde{h} (z (K_{t})^{*}).
\]

In particular, let $T\subset C (n)$ be any finite subset, and consider the sum
\[
S (T) (x) := \sum _{t\in T} h^{*} (t) (x).
\]
Let $z (T)$ be the cycle $\sum _{t\in T} z (t)$.  Then it follows that
\[
S (T) (x) = \tilde{h} (z (T)^{*}) (x).
\]

\subsection{}
Now suppose
\[
\Sigma := \bigcup _{t\in T} t
\]
is \emph{itself} a convex polyhedral cone.  In general, $\Sigma $ will
not be simplicial, but there is nevertheless a well-defined cycle $z
(\Sigma ) \in Z (X)$ associated to $\Sigma $, namely that which is
induced by the polyhedron $K_{\Sigma }$.  The next lemma follows
easily from the previous discussion and the formal properties of our
cycle apparatus.  We omit the simple proof.

\lemma
We have
\[
S (T) (x) = \tilde{h} (z (K_{\Sigma })^{*}) (x).
\]
\endlemma

\section{Admissible units}

\subsection{}
To this point we have not used the fact that $C$ admits an action by a
finite-index subgroup $V$ of $U$, the group of totally positive units.
In this and the next section we use this structure to construct a new fan
$C'$ with $S (C) (x) = S (C') (x)$, and then to construct a sequence of
partial sums for $S (C') (x)$.  Then in the final section we apply the
cycle machinery and Hurwitz's geometric interpretation of $h$ to
complete the proof of Theorem 1.  We begin by discussing some special
collections of units in $V$.  From now on, we assume that $n\geq 3$,
since in the quadratic case a direct proof of the main theorem is easy.

\subsection{}
Recall that $E = (e_{1},\dots ,e_{n})$ is the canonical basis of $\RR
^{n}$.  Let $\varepsilon \in U$ be a totally positive unit.  The
\emph{limit pair} $L (\varepsilon )$ is the pair of projective points
$(\varepsilon (-\infty ), \varepsilon (\infty ))\subset \PP^{n-1}
\times \PP ^{n-1}$, where
\[
\varepsilon (\alpha ) := \lim _{t\rightarrow \alpha } \varepsilon ^{t}.
\]
It is easy to see that $\varepsilon (\alpha )$ always has the form
$\sum _{i\in I } e_{i}$, where we abuse notation slightly and use the
same symbol for a point in $\RR ^{n}$ and the point it induces in $\PP
^{n-1}$.  We say $\varepsilon $ is \emph{generic} if its limit pair
has the form $(e_{i}, e_{j})$, where $i\not =j$.

Let $\sets{n} = \{1,\dots ,n \}$, and let $\cI$ be the set of all
subsets of $\sets{n}$ of order $(n-1)$.  Recall that for any $i\in
\sets{n}$ and $x\in F$, we denote $v_{i} (x)$ by $x^{(i)}$.  In what follows
indexing subscripts and superscripts referring to the real places of
$F$ will be taken modulo $n$.

\definition Let
$T = \{\varepsilon _{1},\dots ,\varepsilon _{n}\}\sub U$ be
a set of totally positive units such that for any $I\in \cI$,
the subset $T_{I} = \{\varepsilon _{i}\mid i\in I\}$
is independent (i.e., the regulator of $T_{I}$ is nonzero).
We say that $T$ is \emph{admissible} if the following hold:
\slist
\lister For each
$\varepsilon \in T$, the coordinates $(\varepsilon ^{(1)},\dots
,\varepsilon ^{(n)})$ are distinct.
\lister We have $L (\varepsilon
_{i}) = (e_{i}, e_{i+1})$.
\lister We have $L (\varepsilon
_{i}/\varepsilon _{j}) = (e_{i}, e_{j})$ for $i\not =j$.
\elist

\lemma
Let $\varepsilon_{1},\dots ,\varepsilon_{n}\in V$ and let $T =
\{\varepsilon_{1}, \dots , \varepsilon_{n} \}$.
Suppose that $T$ satisfies the following: there exist real numbers
$b>a>1$ such that for each $\varepsilon _{i}\in T$, we have
\slist
\lister$\varepsilon _{i}^{(i)} < 1$, $\varepsilon _{i}^{(j)} > 1$ if
$j\not =i$;
\lister $\varepsilon _{i}^{(i+1)} > \varepsilon _{i}^{(i+2)} > \cdots
> \varepsilon _{i}^{(i-1)}$;
\lister $\varepsilon _{i}^{(j)}/\varepsilon _{i}^{(k)} \in (a^{-1}, a)$
for all $j,k \not = i$; and
\lister $\varepsilon _{i}^{(j)}/ \varepsilon _{i}^{(i)} > b$, for all
$j\not =i$.
\elist
Then $T$ is admissible.
\endlemma

\proof The first admissibility condition is clearly satisfied, and we
need only check that limit pairs behave as desired.  The condition $L
(\varepsilon _{i}) = (e_{i}, e_{i+1})$ is obvious. The condition $L
(\varepsilon _{i}/\varepsilon _{j}) = (e_{i}, e_{j})$ follows since in
the ratio $\mu = \varepsilon _{i}/ \varepsilon _{j}$, the smallest
(resp. largest) component is $\mu ^{(i)}$ (resp. $\mu ^{(j)}$).
\endproof

\proposition
Any finite-index $V\subset U$ contains an admissible set of units.
In fact, for every $b>a>1$, there exists $T\sub V$ satisfying
the hypothesis of lemma 3.
\endproposition

\proof  Consider the standard map $\log \colon
\Delta _{\infty }\rightarrow \RR ^{n}$ given by $x \mapsto (\log
x^{(1)},\dots, \log x^{(n)} )$.  The group $V$ is taken to a discrete
subgroup $L (V) \subset \RR ^{n}$, which we may view as a lattice in the
hyperplane $H = L (V)\otimes \R $.  This also endows $H$ with a $\QQ
$-structure, namely $H (\Q ) = L(V) \otimes \Q$.

Now fix $b>a>1$ and consider what conditions the hypotheses of Lemma 3
become in the subspace $H$.  Let $(\xi _{1}, \dots , \xi _{n})$ be the
coordinates of a point in $\R ^{n}$ with respect to $E$, so that $H$
is defined by the equation $\sum \xi _{i} = 0$.  (Note that the
rational structure induced by $E$ is \emph{not} the same
as that induced by $L (V)$.)  We find that a set of units
$\{\varepsilon _{1},\dots ,\varepsilon _{n} \}$ satisfies (1)--(4) of
the statement if and only if for each $i$, the point $(\log
\varepsilon _{i}^{(1)}, \dots , \log \varepsilon _{i}^{(n)})$ lies in
the region $R_{i}\subset H$ determined by the inequalities
\slist
\lister $\xi_{i} < 0$, $\xi_{j} > 0$ for $j\not =i$;
\lister $\xi _{i+1} > \xi _{i+2} > \cdots > \xi _{i-1}$;
\lister $-\log a < \xi _{j}- \xi _{k} < \log a$,
for all $j,k \not = i$; and
\lister $\xi _{j} - \xi _{i}> \log b$, for all $j\not =i$.
\elist
We claim that each $R_{i}$ is an unbounded open set in $H$ of full
dimension.  This can be seen as follows.  For fixed $i$, the
conditions (1) define an
$(n-1)$-dimensional simplicial cone
\[
\sigma_{i} = \bigl\{ \sum_{i\not =j} \lambda_{j} (e_{j}-e_{i})\bigm |
\lambda_{j}\in \RR_{\geq 0}\bigr \}.
\]
The
conditions (2) cut out an $(n-1)$-cone $\sigma _{i}'$ in the
barycentric subdivision of $\sigma_{i} $.  Let $\tau_{i} $ be the
barycenter of $\sigma_{i} $, i.e. the $1$-cone
\[
\tau_{i} = \R_{\geq 0} (- (n-1)e_{i} + \sum_{j\not =i} e_{j}),
\]
and let $U_{i}$ be a small
tubular neighborhood of $\tau_{i} $.  Then each inequality in (3)
determines a half-space containing $U_{i}$.  Hence the
intersection of these half-spaces with $\sigma' _{i}$ is unbounded.
Finally, the half-spaces determined by (4) include all points in
$U_i $ that are sufficiently far away from the origin, and hence the
region $R_{i}$ is unbounded and has dimension $n-1$.

Figure 1 shows the $R_{i}$ for the case $n=4$. The vertices of the
cuboctahedron $P$ are the 12 points $\{e_{i}-e_{j} \mid i\not =j\}$.
Each $\sigma _{i}$ is a cone generated by a triangular face of $P$
with apex at the center of $P$.  The $\tau_{i}$ are the rays through
the centers of these four faces of $P$, and the $R_{i}$ are the four
thin semiinfinite prisms.

To conclude the proof, we claim that each $R_{i}$ contains infinitely
many points of $L (V)$.  Indeed, we assert that for each $i$, we have
$L (V)\cap \tau_{i} = \{0 \}$.  This proves the claim, since it
implies the image of $L (V)$ is dense in the quotient $H/\R
\tau_{i}$, and hence the inverse image of any open set in this
quotient must contain infinitely many points of $L (V)$.

To prove the assertion, let $x\in L (V)\cap \tau_{i}$ be
nonzero.  Then we may write
\[
x = (\lambda , \dots ,\lambda , - (n-1)\lambda , \lambda , \dots ,\lambda),
\]
for some positive real number $\lambda$, where the $- (n-1)\lambda $ appears
in the $i$th position.  This implies that there is a unit
$\varepsilon$ of the form
\[
\varepsilon = (e^{\lambda}, \dots , e^{\lambda}, e^{- (n-1)\lambda},
e^{\lambda}, \dots ,
e^{\lambda});
\]
in other words, all $v_{j}(\varepsilon)$ with $j\not = i$ are equal.
We claim that no element of $F$ that is not in $\QQ$ can have this
form.  To see this, let $\alpha\not =\beta$ be the distinct infinite
places of $\varepsilon$, and let $f (X) = (X-\alpha)^{n-1}
(X-\beta)\in \RR [X]$.  Clearly $f (\varepsilon ) = 0$, and so the
minimal polynomial $g (X)$ of $\varepsilon$ must divide $f$.  This
implies $g (X) = (X-\alpha) (X-\beta)$, which implies $\varepsilon\in
F$ is quadratic over $\QQ$.  But this means that under $F\rightarrow \RR^{n}$,
half of the embeddings must equal $\alpha$ and half must equal
$\beta$.  This contradicts the
assumption that $n\geq 3$ and the proof of the proposition
is complete.
\endproof

\vbox{
\epsfysize=2.0in
\centerline{\epsfbox{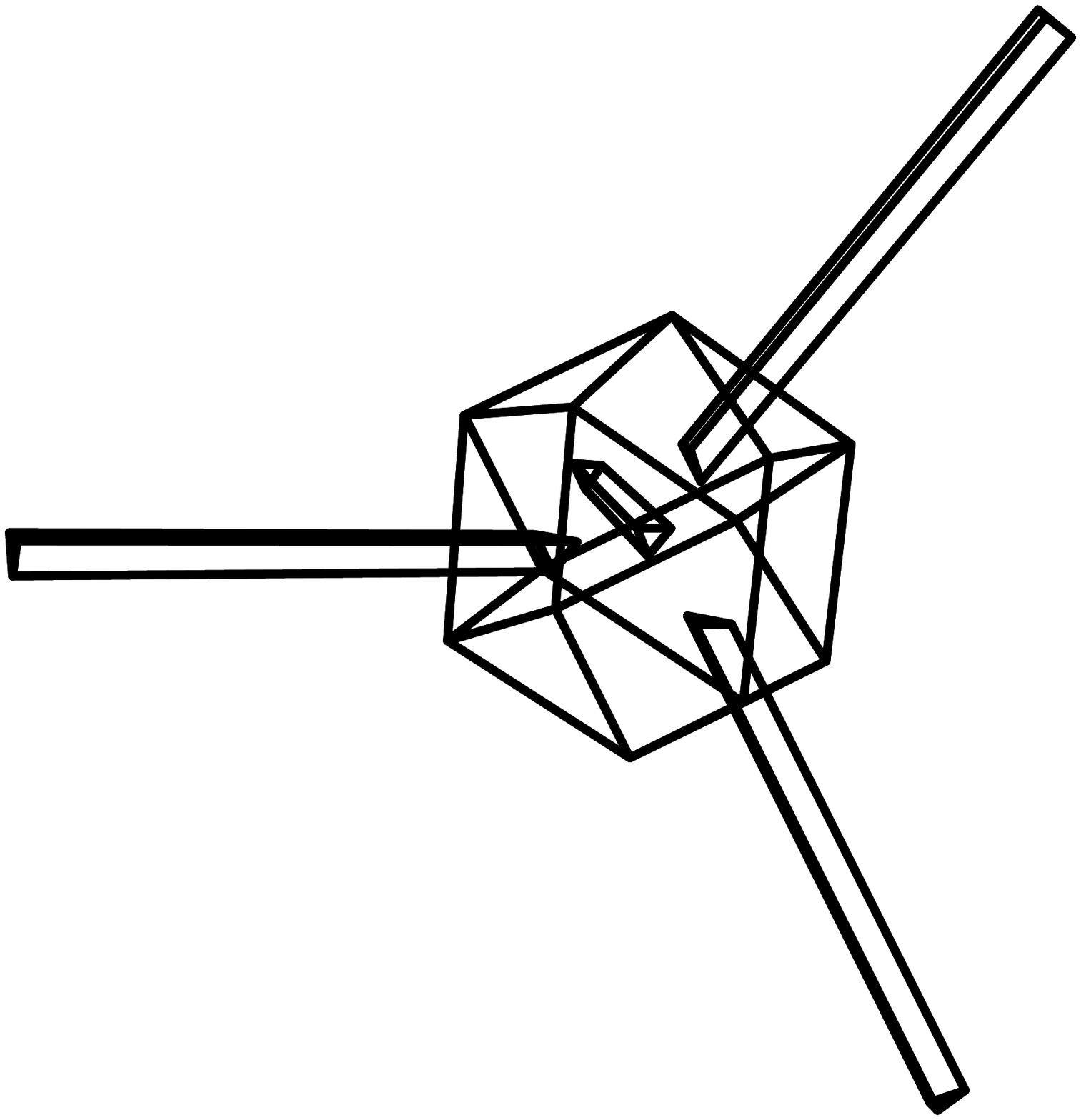}}
\centerline{Figure 1.  The regions $R_{i}$ for a quartic field.}
\bigskip
}
\section{Partial sums}

\subsection{}
Let $V$ be a fixed subgroup of finite index in $U$, the
group of totally positive units.
Choose and fix an admissible set of units $T = \{\varepsilon _{1},
\dots , \varepsilon _{n} \}$ satisfying the conditions of Lemma 3,
with $a,b$ chosen so that $b>a^n>1$.
For each $I\in \cI $, define $ V_I$ by
\[
V_{I} = \bigl\{ \prod_{i\in I} \varepsilon _{i}^{\alpha _{i}}\bigm |
\alpha _{i}\in \ZZ , \sum_{i\in  I} \alpha _{i} = 0\bigr\}.
\]
Let
$\Sigma _I \subset \Delta _{\infty
}$
be the convex cone generated by the half-lines through the points
$
V_{I}.
$

\lemma
The points in $V_I$ generate the spanning rays of $\Sigma
_I$.  In other words, no point of $V_I$ lies in the
relative interior of $\Sigma_I$.
\endlemma

\proof
Let $j = \sets{n} \smallsetminus I$ and
let $\varphi_j\colon  \D_\i\ra \R^{n-1}_{+} $
be the map
$$ \varphi_j(x_1,\dots ,x_n) \ = \ \left(
{x_1\over x_j}, \dots ,\widehat{\frac{x_{j}}{x_{j}}} , \dots ,{x_n\over x_j}\right)
$$
where the $j^{th}$ component is omitted.
Since $\varphi_j$ takes
straight lines to straight lines, we see that
the image
$\Pi_I = \varphi_j (
\Sigma_I )$ is a convex polyhedral
subset of $\RR_+ ^{n-1}$ with $\varphi_j^{-1}(\Pi_j)=\Sigma_I$.  To
prove the lemma it thus suffices to show that the points
\[
\bigl\{\varphi_j (v) \bigm| v\in V_{I}\bigr\}\eqno (\showeqnum)
\]
are the vertices of $\Pi = \Pi_I $.

Now consider the set $\partial \Pi _{\RR}$ of the boundary of the
``real points'' of $\Pi $, by which we mean the subset in $\RR ^{n-1}$
of points of the form
$$
\partial \Pi _{\RR} := \bigl\{ \varphi_{j} \left(\prod_{i\in I} \varepsilon
_{i}^{t _{i}}\right)\bigm |t _{i}\in \RR ,
\ \ \sum_{i\in I} t _{i} =
0\bigr \}.
$$
We claim that
$$
\pl\Pi _{\RR } = \bigl\{(z_{1},\dots ,z_{n-1}) \bigm | z_{i} > 0, \prod
z_{i}^{a_{i}} = 1 \bigr \},
$$
where $(a_{1},\dots ,a_{n-1})$ are certain fixed real
numbers.  To see this, assume (to ease notation) that $j=n$, and let
$Z$ (respectively $T$) be the column vector $\ ^t(z_1,\dots,z_{n-1})$ (resp. $\ ^t(t_1,\dots,t_{n-1})$).
Then $Z\in \pl\Pi_\R$ if and only
if $\log Z = E \cdot T$ for some $T$, where $E$ is the matrix whose
$p,q$ entry is 
\[
(\log \e^{(p)}_q - \log \e^{(j)}_q), \quad p,q=1,\dots ,n-1.
\]
Since the $\e_i$ generate a subgroup of the unit
group with full rank, the matrix $E$ is invertible.  This
proves the existence of the $a_i$. In fact,
we have $(a_1,\dots, a_{n-1}) = (1,\dots,1)E^{-1}$.  The same argument
works for any $j=1,\dots ,n-1$, which proves the claim.

Next we claim that all the $a_i$ are positive: If not,
then we can rewrite $\p z_i^{a_i}=1$ as
$\p_{k\in A} z_k^{r_k} \ = \ \p_{l\in B} z_l^{s_l} $
where $I=A\cup B$ is  a partition into disjoint non-empty
subsets, and the $r_k, s_l$ are all positive real numbers.
Let $r_{k_0}=\max r_k$ and $s_{l_0}=\max s_l$, and let
$z = \varphi(\e_{l_0}/\e_{k_0})$. Then $z\in \pl\Pi_\R$. Moreover,
$z_{k_0} > (ba^{-1})^{r_{k_0}}$ and $z_k > a^{-r_{k_0}}$ if
$k\not=k_0$. Since $b>a^n$, this shows 
$\p_{k\in A} z_k^{r_k} >1$. Similarly, $\p_{l\in B} z_l^{s_l}<1$.
This is a contradiction, and thus we conclude that
all the $a_i$ are positive.

Let $\Pi _{\RR }$  be defined as
$$
\Pi _{\RR } = \bigl\{(z_{1},\dots ,z_{n-1}) \bigm | z_{i} > 0, \prod
z_{i}^{a_{i}} \geq 1 \bigr \},
$$
We claim that $\Pi _{\RR }$ is a convex subset of $\RR ^{n-1}$.
Note first that
$\partial \Pi _{\RR }$ is the graph of the function
$f\colon \RR ^{n-2} \rightarrow \RR $ defined by
$$
f (z_{1},\dots ,z_{n-2}) = \prod _{i=1}^{n-2} z_{i}^{-p_{i}}, \quad p_{i}>0.
$$
Now consider the Hessian matrix
\[
M = \Bigl(\frac{\partial^{2} f}{\partial z_{i}\partial z_{j}}\Bigr).
\]
We claim that $M$ is positive-definite, from which the convexity of
$\Pi _{\RR }$ follows immediately.  To see this, for $k=1,\dots ,n-2$
let $M_{k}$ be the $k\times k$ upper-left submatrix of $M$.  Then by
an induction argument
$$
\det M_{k} = (1+\sum_{i\leq k} p_{i} )\prod_{i\leq k} p_{i}z_{i}^{- (n-2) p_{i}-2}.\eqno (\showeqnum)
$$
Since $p_{i}, z_{i}>0$, this determinant is positive on $\partial \Pi
_{\RR }$ for each $k$.  By a standard result of linear algebra, this
implies $M$ is positive-definite, and thus $\Pi _{\RR }$ is a convex
subset of $\RR ^{n-1}$.

We claim convexity of $\Pi _{\RR } $ implies that the points in
\eqref{4} are exactly the vertices of $\Pi $.  To see this, assume
that $\varphi_j(v)$ is in the relative interior of $\Pi$
for some $v\in V_I$. Since $\Pi _{\RR }$
is convex, we deduce that $\Pi\sub\Pi_\R$ and hence
$\varphi_j(v)$ is in the relative interior of $\Pi_\R$.
But $\varphi_j(v)\in \pl\Pi_\R$.
This
contradiction completes the proof.  \endproof

\lemma

\slist 
\lister The faces on the $n-1$ dimensional convex polyhedron
$\Pi_I\sub \R^{n-1}_+$
are compact polyhedra of dimension $n-2$.

\lister If
we let $\cF_I$ denote the set of faces of $\Pi_I$, then there
exists a finite subset $S_I\sub \cF_I$ with the following property:
For every $f\in \cF_I$  there exists $\e\in V_I$ such
that $f = \e\cdot s$ for some $s\in S_I$.
\elist 
\endlemma

\proof
Let $f\in \cF$, let $H$ be the hyperplane spanned by $f$, and $v\in
f\sub H$ be a vertex of~$f$.  For each $i = 1,\dots ,n$, choose a
$\e_i\in V_I$ that is very close to the standard basis element $e_i$
inside $\P^{n-1}$ (in other words, the $i^{th}$ component of $\e_i$ is
much larger than all the other components). For $r$ a positive
integer, we let $K(r)$ be the hyperplane spanned by the $\e_i^r$.  We
see that $K(r)\cap \R^{n-1}_{\geq 0}$ is a compact simplex of
dimension $n-2$.

Let $\tilde K(r)$ be the $n-1$ simplex spanned by $K(r)$ and the
origin. Then $\tilde K(r)$ is compact and contains $H\cap
\R^{n-1}_{\geq 0}$ for $r$ sufficiently large. Thus shows that $H\cap
\R^{n-1}_{\geq 0}$ is a compact simplex, and thus $H\cap \Pi$ is
compact, which proves the first part of Lemma 5.

Now we note that $V_I$ acts on $\pl\Pi_\R$ and that $\pl\Pi_\R/V_I$ is
compact.  Thus we can choose $r$ such that $\tilde K(r)$ contains a
compact fundamental domain $D\sub \pl\Pi_\R$.  Now let $f\in \cF_I$,
let $v\in f$ be a vertex, and choose $\e\in V_I$ such that $\e\cdot  v\in
D$. Then $\e \cdot f $ is a face of $\pl\Pi$, and by the proof of part one,
$\e\cdot f\sub \tilde K(r')$ for some $r'>r$ that depends on $D$ but not
on $f$. Thus all the vertices of $\e\cdot f$ lie in $\tilde K(r')\cap
\Pi_R$, which is compact.  Since $V_I\cap \tilde K(r')\cap \Pi_R$ is
finite, the second part of the lemma is proved.
\endproof

\subsection{}
Now we come to the main tool that we need to prove Theorem 1.

\theorem
Let $C$ be a good fan, and $V \sub U$ a subgroup of finite
index acting on $C$.
Then
there
exists $C'$, a good refinement of $C$, and a family of
convex subsets $\Sigma_N\sub \D_\i$, with the following
properties:

1. For each $N$, the set $\Sigma_N$ is a finite union of
top dimensional simplices in $C'$.

2. $\D_\i = \bigcup_{_N}\Sigma_N$.

\endtheorem

\proof
 Let $\{\e_1,\dots ,\e_n\}$ satisfy the
conditions of Lemma 3, with $a,b$ chosen so that
$b>a^n>1$.
For each $j, 1\leq j\leq n$, let
$I(j)=\{i\mid 1\leq i\leq n, i\not=j\}$.
Define
$$ \Sigma_N = \bigcap_{J=1}^N \e_j^{-N}\cdot\varphi_j^{-1}(\Pi_{I(j)})
$$

Now fix $j$, and let $I=I(j)$. Let $\cD_j=\varphi_j(C)$, which is a
simplicial  decomposition of $\R^{n-1}$. Choose a fundamental domain
$F$ for the action of $V$ on $\R^{n-1}$ of the form
$$ F = \bigcup_{t\in T} t,
$$
where $T$ is a finite subset of top dimensional simplices
in $\cD$, and
let
$S_I$ be as in Lemma~5. Since $S_I$ is finite, and since, by
Lemma 5, the elements of $S_I$ are compact, the set
$$ \bigl\{(\e, t, s)\in V\times T\times S_{I} \bigm | \e\cdot s\cap t \not=
\emptyset\bigr\}
$$
is finite. Thus there is a simplicial decomposition $\cD_j'$,
which is a refinement of $\cD_j$, with the property that every
$f\in \cF_I$ is a finite union of $n-2$ simplices in $\cD_j'$.
\v
Let $C'$ be a common refinement of $\varphi_J^{-1}(\cD_j')$
 for all $j$. Then $C'$ clearly satisfies condition 1 of
Theorem 3, and $\Sigma_N$ clearly satisfies condition 2.
\endproof

\section{Proof of the main theorem}

We now complete the proof of Theorem 1.  Let $x_{0}\in F$ be totally
positive, and consider the sequence
of cones $\{\Sigma _{N} \}$ constructed in Theorem 3.  By Theorem 3 we have $\lim
_{N\rightarrow \infty } S (\Sigma _{N}) (x_{0}) = S (C') (x_{0})$,
which in turn equals $S (C) (x_{0})$ by Lemma 1.  We want to show
\[
\lim _{N\rightarrow \infty } S (\Sigma _{N}) (x_{0})  = 1/\Norm (x_{0}).
\]
Since the sets $\Sigma _{N}$ exhaust $\Delta _{\infty }$, we have
$x_{0}\in \Sigma _{N}$ for sufficiently large $N$.  Because $\Sigma _{N}$ is
convex, we know that the singular hyperplanes of the rational function
$S (\Sigma _{N}) (x)$ miss $x_{0}$.  Hence the sequence of partial
sums in the limit is well-defined.

The cycle machinery from \S\S 3--7 implies
\[
S (\Sigma _{N}) (x_{0}) = \tilde{h} (z (\Sigma _{N})^{*}) (x_{0})
\]
and
\[
z (\Sigma _{N})^{*} = z (\Sigma _{N}^{*}).
\]
Hence we can compute $S (\Sigma _{N}) (x_{0})$ by arbitrarily dividing $\Sigma
_{N}^{*}$ into \emph{simplicial} cones
$$
\Sigma _{N}^{*} = \bigcup
_{\sigma \in T} \sigma,
$$
where $T$ is some finite set of top dimensional cones (not in $C$, of
course).  We can then apply the relation in $Z (\PP ^{n-1})$
\[
z (\Sigma ^{*}_{N}) = \sum _{\sigma \in T} z (\sigma ).
\]

Now we use Hurwitz's interpretation of the function $h$.  For any
cone $\sigma \subset \RR ^{n}$, let $\PP \sigma $ be the induced set
in $\PP ^{n-1}$.  Let $H_{x_{0}}$ be the hyperplane in $\PP^{n-1}$
determined by the hyperplane in $\RR^{n}$ orthogonal to $x_{0}$.  If
$N$ is sufficiently large, then $x_{0}\in \Delta_{N}^{*} $ and $H_{x_{0}}$ misses $\PP
\Delta _{N}^{*}$.  Thus $H_{x_{0}}$ eventually misses $\PP \Sigma
_{N}^{*}$ and the sets $\{ \PP \sigma \mid \sigma \in T\}$.
This implies that all three sets
$\PP \Sigma _{N}^{*}$, $\PP \Delta _{N}^{*}$, and $\PP \Delta _{\infty
}$ are contained in the affine chart $\PP ^{n-1} \setminus H_{x_{0}}$.
Therefore we may fix an orientation such that
\[
S (\Sigma _{N}^{*}) (x_{0}) = \int _{\PP \Sigma _{N}^{*}} \Omega _{x_{0}},
\]
and similarly for $\Delta _{\infty }$ and $\Delta _{N}^{*}$.
Since $\Delta _{N}^{*} \supset \Sigma _{N}^{*} \supset \Delta _{\infty
}$ and $\Delta _{N}^{*}\rightarrow \Delta _{\infty }$ as $N\rightarrow
\infty $ (Theorem 3), we have
\[
\int _{\PP \Sigma _{N}^{*}} \Omega _{x_{0}} \longrightarrow
\int _{\PP \Delta _{\infty }} \Omega _{x_{0}}.
\]
Since $\int _{\PP \Delta _{\infty }} \Omega _{x_{0}} = 1/\Norm (x_{0})$,
Theorem 1 is proved.

\section{References}

\def\bibitem#1&{\hangindent4em\hangafter=1\hbox to \hangindent{{[#1]}\hss}\ignorespaces}
\everypar={\bibitem}

AMRT&A.~Ash, D.~Mumford, M.~Rapoport, and Y.~Tai, \emph{Smooth
compactifications of locally symmetric varieties}. Math. Sci. Press,
Brookline, Mass., 1975.

ADS&M. Atiyah, H. Donnelley, and I. M. Singer, \emph{Eta invariants,
signature defects of cusps, and values of ${L}$-functions}. Ann. of
Math. {\bf 118} (1983), pp. 131--177.

Ful&W.~Fulton, \emph{Introduction to toric varieties}. Ann. of
Math. Studies 131, Princeton 1993.

GP&H.~Cohen et~al., \emph{The {GP}-{P}ari system}, available from
  {\tt ftp://megrez.math.u-bordeaux.fr/}.

Hir&F.~Hirzebruch, \emph{Hilbert modular surfaces}.  Enseignement
Math. (2) {\bf 19} (1973), 183--281.

Hur&A.~Hurwitz, \emph{\"Uber die Anzahl der Klassen positiver
tern\"arer quadratischer Formen von gegebener Determinante}.
Math. Ann. {\bf 88} (1923), pp. 26--52.

KKMS&G.~Kempf, F.~Knudsen, D.~Mumford, and B.~St. Donat,
\emph{Toroidal embeddings I}.  Lecture notes in math. {\bf 339}
(1973), Springer-Verlag.

Mey&C.~Meyer, \emph{Die Berechnung der Klassenzahl abelischer K\"orper
\"uber quadratischen Zahlk\"orpern}. Berlin 1957.

M\"ul&W.~M\"uller.  \emph{Manifolds with cusps of rank one}.  Lecture
notes in math. {\bf 1244} (1987), Springer-Verlag.

Oga&S.~Ogata, \emph{Generalized Hirzebruch's conjecture for
Hilbert-Picard modular cusps}.  Japan. J. Math. (N.S.)  {\bf 22} (1996),
no. 2, 385--410.


Sa&I.~Satake, \emph{On zeta functions associated with cones and cusp
singularities of the second kind}, in Automorphic forms and related
topics, Jin-Woo Son and Jae-Hyun Yang, eds., pp. 17--26, Pyungsan
Institute for Mathematical Sciences, Seoul, Korea.

Sc1&R.~Sczech, \emph{Eisenstein group cocycles for $GL_{n}$ and values
of $L$-functions}.  Invent. math. {\bf 113} (1993), pp. 581--616.

Sc2&R.~Sczech, \emph{Intersection numbers of cusp divisors
and values of Hecke $L$-functions}.  Preprint 2001.

Shi&J.-Y. Shi, \emph{The Kazhdan-Lusztig cells in certain affine Weyl
groups}.  Lecture notes in math. {\bf 1179}, Springer-Verlag.

Zag&D.~Zagier, \emph{Nombres de classes et fractions
continues}. Journ\'ees Arithm\'etiques de Bordeaux (Conf.,
Univ. Bordeaux, Bordeaux, 1974), pp. 81--97. Asterisque, No. 24-25.

\everypar={}
\bigskip
\noindent
\line{Department of Mathematics and Statistics\hfill }
\line{University of Massachusetts\hfill }
\line{Amherst, MA 01003\hfill }
\line{\emph{gunnells@math.umass.edu}}
\bigskip
\noindent
\line{Department of Mathematics and Computer Science\hfill }
\line{Rutgers University\hfill }
\line{Newark, NJ 07102--1811\hfill }
\line{\emph{sturm@andromeda.rutgers.edu}}

\end